\documentclass[a4paper,12pt]{amsart}
\usepackage{amssymb}
\usepackage{ifthen}
\usepackage{graphicx}
\usepackage{float}
\usepackage{caption}
\usepackage{subcaption}
\usepackage{cite}
\usepackage{amsfonts}
\usepackage{amscd}
\usepackage{amsxtra}
\usepackage{mathrsfs}
\usepackage{mathrsfs}
\usepackage[usenames]{color}

\newtheorem{thm}{Theorem}
\newtheorem{cor}{Corollary}
\newtheorem{lem}{Lemma}

\newtheorem{rem}{Remark}
\newtheorem{example}{Example}
\newtheorem{defn}{Definition}
\newtheorem{prob}{Problem}

\newtheorem{conj}{Conjecture}
\theoremstyle{definition}

\newcounter {own}
\def\theown {\thesection  .\arabic{own}}

\newenvironment{pf}[1][]{%
 \vskip 3mm
 \noindent
 \ifthenelse{\equal{#1}{}}%
  {{\slshape Proof. }}%
  {{\slshape #1.} }%
 }%
{\qed\bigskip}

\newcounter{alphabet}
\newcounter{tmp}

\makeatletter
\newcommand{\Ref}[1]{\@ifundefined{r@#1}{}{\setcounter{tmp}{\ref{#1}}\Alph{tmp}}}
\makeatother

\newenvironment{Lem}[1][]{\refstepcounter{alphabet}%
\bigskip%
\noindent%
{\bf Lemma \Alph{alphabet}}%
{\bf .} \itshape}{\vskip 8pt}

\newcommand{\IR}{{\mathbb R}}
\newcommand{\IN}{{\mathbb N}}
\newcommand{\IC}{{\mathbb C}}
\newcommand{\ID}{{\mathbb D}}

\newcommand{\D}{{\mathbb D}}

\def\be{\begin{equation}}
\def\ee{\end{equation}}

\newcommand{\bee}{\begin{enumerate}}
\newcommand{\eee}{\end{enumerate}}

\newcommand{\blem}{\begin{lem}}
\newcommand{\elem}{\end{lem}}
\newcommand{\bthm}{\begin{thm}}
\newcommand{\ethm}{\end{thm}}
\newcommand{\bcor}{\begin{cor}}
\newcommand{\ecor}{\end{cor}}
\newcommand{\beg}{\begin{example}}
\newcommand{\eeg}{\end{example}}
\newcommand{\begs}{\begin{examples}}
\newcommand{\eegs}{\end{examples}}
\newcommand{\bdefe}{\begin{defn}}
\newcommand{\edefe}{\end{defn}}
\newcommand{\bprob}{\begin{prob}}
\newcommand{\eprob}{\end{prob}}
\newcommand{\bei}{\begin{itemize}}
\newcommand{\eei}{\end{itemize}}

\newcommand{\bcon}{\begin{conj}}
\newcommand{\econ}{\end{conj}}
\newcommand{\bcons}{\begin{conjs}}
\newcommand{\econs}{\end{conjs}}
\newcommand{\bprop}{\begin{propo}}
\newcommand{\eprop}{\end{propo}}
\newcommand{\br}{\begin{rem}}
\newcommand{\er}{\end{rem}}
\newcommand{\brs}{\begin{rems}}
\newcommand{\ers}{\end{rems}}
\newcommand{\bo}{\begin{obser}}
\newcommand{\eo}{\end{obser}}
\newcommand{\bos}{\begin{obsers}}
\newcommand{\eos}{\end{obsers}}
\newcommand{\bpf}{\begin{pf}}
\newcommand{\epf}{\end{pf}}
\newcommand{\ba}{\begin{array}}
\newcommand{\ea}{\end{array}}
\newcommand{\beq}{\begin{eqnarray}}
\newcommand{\beqq}{\begin{eqnarray*}}
\newcommand{\eeq}{\end{eqnarray}}
\newcommand{\eeqq}{\end{eqnarray*}}

\newcommand{\ds}{\displaystyle}

\newcounter{minutes}
\divide\time by 60
\newcounter{hours}
\multiply\time by 60 \addtocounter{minutes}{-\time}

\begin{document}
\bibliographystyle{amsplain}
\title[Coefficients of univalent harmonic mappings]{Coefficients of univalent
harmonic mappings}

\thanks{
File:~\jobname .tex,
          printed: \number\year-\number\month-\number\day,
          \thehours.\ifnum\theminutes<10{0}\fi\theminutes}

\author{Saminathan Ponnusamy
}
\address{S. Ponnusamy, and A. Sairam Kaliraj,
Indian Statistical Institute (ISI), Chennai Centre, SETS, MGR Knowledge City, CIT
Campus, Taramani, Chennai 600 113, India. }
\email{samy@iitm.ac.in, samy@isichennai.res.in; sairamkaliraj@gmail.com}
\author{Anbareeswaran Sairam Kaliraj}

\author{Victor V. Starkov}
\address{V.V. Starkov, Petrozavodsk State University, 33, Lenin Str., 185910, Petrozavodsk, Republic of Karelia, Russia.}
\email{vstarv@list.ru}


\date{\today  
}

\begin{abstract}
Let $\mathcal{S}_H^0$ denote the class of all functions $f(z)=h(z)+\overline{g(z)}=z+\sum^\infty_{n=2} a_nz^n +\overline{\sum^\infty_{n=2} b_nz^n}$ that are sense-preserving, harmonic and univalent in the open unit disk $|z|<1$. The coefficient conjecture for $\mathcal{S}_H^0$ is still \emph{open} even for $|a_2|$. The aim of this paper is to show that if $f=h+\overline{g} \in \mathcal{S}^0_H$ then $ |a_n| < 5.24 \times 10^{-6} n^{17}$ and $|b_n| < 2.32 \times 10^{-7}n^{17}$ for all $n \geq 3$. Making use of these coefficient estimates, we also obtain radius of univalence of sections of univalent harmonic mappings.
\end{abstract}
\thanks{ }

\maketitle
\pagestyle{myheadings}
\markboth{S. Ponnusamy, A. Sairam Kaliraj, and V. V. Starkov}{Coefficients of univalent harmonic mappings}

\section{Preliminaries and main results}\label{PSS8Sec1}
Let $\ID$ denote the open unit disk centered at the origin of the complex plane. The theory of univalent and sense-preserving
complex-valued harmonic functions in $\ID$ has attracted a lot of attention since the appearance of the paper by  Clunie and Sheil-Small
\cite{Clunie-Small-84} which brought the theory a large step forward. They pointed out that many of the classical results for
conformal mappings have clear analogues for harmonic mappings although only few of them have been addressed and used by a number of authors, while
others were not because of the higher difficulty level. Besides its interest from the point of view of analysis, it
has been recently shown to be of relevance in some problems related to fluid flows. This
applied mathematics connection brings new relevance to the issue of coefficient estimates
for a family of sense-preserving harmonic mappings since these maps provide an approach towards
obtaining explicit solutions to the incompressible two-dimensional Euler equations.
For example, A.~Aleman and A.~Constantin \cite{AleConst2012} pointed out the importance of harmonic mappings in the Eulerian
description of fluid flows and developed a method which is largely based
on a detailed study of the governing equations using analytic function theory, and an important role played
by the univalence of the labelling map.   Also, the authors in \cite{AleConst2012} presented several examples to illustrate
how the classical solutions can be obtained from the more general solution formulas via univalent
harmonic mappings. More recently, O.~Constantin and M.J.~Mart\'{i}n \cite{ConstMartin2017} continued this investigation and
proposed a different approach that provides a complete
solution to the original problem of classifying all two-dimensional ideal fluid flows with harmonic Lagrangian labelling mappings.
This approach is based on the ideas from the theory of planar harmonic mappings and thus, provide an illustration of the
deep link between the sense-preserving harmonic mappings and fluid flow problems. This
newly explored connection renewed our interest in this topic.

In this article, we consider the class $\mathcal{S}_H$ of all univalent, sense-preserving harmonic functions $f$ of $\ID$ normalized by
$f(0)=0=f_z(0)-1$. Every such function has a unique canonical representation of the form $f=h+\overline{g}$, with $h$ and $g$ analytic in $\ID$ and $g(0)=0$.
Here $h$ and $g$ are often referred to as analytic and co-analytic parts of $f$.
Let $\mathcal{S}^0_{H}  = \{f=h+\overline{g} \in \mathcal{S}_{H}:\, g'(0)=0 \}$. Clearly,
$\mathcal{S}= \{f=h+\overline{g} \in \mathcal{S}_{H}:\, g(z)\equiv 0 \}$ is the class of normalized univalent analytic functions in $\ID$.
A typical element $f \in \mathcal{S}^0_H$ has the form
\be\label{PSSerRep}
f(z)=h(z)+\overline{g(z)}:=z+\sum _{k=2}^{\infty}a_k z^k + \sum _{k=2}^{\infty}\overline{b_k z^k}, ~z \in \mathbb{D}.
\ee
Throughout the discussion we shall use this representation. Clunie and Sheil-Small \cite{Clunie-Small-84} proved that both $\mathcal{S}_{H}$ and
$\mathcal{S}^0_{H}$ are normal whereas only $\mathcal{S}^0_{H}$ is compact with respect to the topology of uniform convergence on compact subsets of $\D$.
A function $f \in \mathcal{S}^0_H$ is said to belong to the class $\mathcal{S}^{*0}_{H}$, $\mathcal{K}^0_{H}$ and
$\mathcal{C}^0_{H}$ if $f(\mathbb{D})$ is starlike with respect to the origin, convex and close-to-convex, respectively.
The corresponding notations for the analytic case are $\mathcal{S}^*$, $\mathcal{K}$ and $\mathcal{C}$, respectively.
For basic information about $\mathcal{S}^0_{H}$ and related geometric subfamilies,
one can refer to \cite{Clunie-Small-84,Duren,Duren:Harmonic} and the recent expository article of Ponnusamy and Rasila \cite{PonRasi2013}.

This article is organized as follows. In Section \ref{PSS8SubSec1}, we present a preliminary information on the coefficient conjecture of
Clunie and Sheil-Small \cite{Clunie-Small-84} and present a coefficient estimate for a family of sense-preserving harmonic mappings, which contains the class $\mathcal{S}^0_H$ (Theorem \ref{PS8Thm1}). In Section \ref{PSS8SubSec2}, we recall some known results on the sections of functions in certain geometric subclasses of univalent harmonic mappings and present our result on the radius of univalence of partial sums of functions in $\mathcal{S}^0_H$ (Theorem \ref{PS8Thm2}). Few basic lemmas that are needed for the proofs of these two results are recalled in Section \ref{PSS8Sec2}. The proofs of our main results are presented in Section \ref{PSS8Sec3}. Some consequences of them are discussed in Section \ref{PSS8Sec4} (Theorems \ref{spec_coeff_bound} and \ref{part_spec_case}).

\subsection{Coefficient conjecture of Clunie and Sheil-Small}\label{PSS8SubSec1}
Using the method of shearing, Clunie and Sheil-Small \cite{Clunie-Small-84} obtained an important member of so-called slit mapping $K=H+\overline{G}$,
where
$$H(z)=\frac{z-\frac{1}{2}z^2+\frac{1}{6}z^3}{(1-z)^3}= z+\sum_{n=2}^{\infty}A_nz^n
~\mbox{ and }~G(z)=\frac{\frac{1}{2}z^2+\frac{1}{6}z^3}{(1-z)^3}=\sum_{n=2} ^{\infty}B_nz^n
$$
with
$$A_n=\frac{(n+1)(2n+1)}{6} ~\mbox{ and }~ B_n=\frac{(n-1)(2n-1)}{6} ~\mbox{ for $n\geq 2$.}
$$
The function $K$  is called the harmonic Koebe function and it maps the unit disk one-to-one onto the slit domain
$\IC\backslash\{u+iv:\, u\leq -1/6,\,v=0\}$ which is indeed convex along horizontal direction,
and it plays an extremal role for several extremal problems in $\mathcal{S}^{*0}_{H}$ and $\mathcal{C}_{H}^0$, such as coefficient bounds and
covering theorems (see \cite{Clunie-Small-84, Duren:Harmonic, Sheil-Small}). Due to the extremal role of the harmonic Koebe function in these
families, it was natural for Clunie and Sheil-Small \cite{Clunie-Small-84} to conjecture that if
$f=h+\overline{g}\in {\mathcal S}_{H}^{0}$ is given by \eqref{PSSerRep},
then for all $n\geq2$,
$$|a_n|\leq A_n,~ |b_n|\leq  B_n~ \mbox{ and }~\big||a_n|-|b_n|\big|\leq  n
$$
and equality occurs for $f(z)=K(z)$. In \cite{Clunie-Small-84}, they also showed that
$|b_2|\leq 1/2$ which is sharp, and the non-sharp estimate $|a_2|<12172$.  Later in 1990, Sheil-Small \cite{Sheil-Small} improved it
to $|a_2|<57$, and then  Duren \cite[p.~96]{Duren:Harmonic} improved it further to $|a_2|<49$ which is again far from the conjectured bound
$|a_2|\leq 5/2$. The above conjecture remains \emph{open} and little is known for $n\geq 3$ for the class ${\mathcal S}_{H}^{0}$. In \cite{Starkov}, it has been proved that
\be\label{coeff_bound_starkov}
|a_n| <  \frac {(2e^2)^\alpha}{2\alpha} n^\alpha ~\mbox{ and }~ |b_n| <  \frac {(2e^2)^\alpha}{2\alpha} n^\alpha ~\mbox{ for all }~ n\in \IN,
\ee
where $\alpha: = {\rm ord} \,\mathcal{S}_H = \sup_{f\in \mathcal{S}_H} |a_2|$. However, finding the explicit value of  $\alpha$ or even finding a good upper
bound itself seems to be a difficult task. Very recently,  Abu Muhanna et. al \cite{AAP-PP2015} obtained the following result, which is the best known bound so far and this could be used in \eqref{coeff_bound_starkov}.

\begin{Lem} {\rm\cite{AAP-PP2015}}\label{samy_abu_bound}
If $f=h+\overline{g} \in \mathcal{S}^0_H$, then $|a_2|=|h''(0)/2|\leq 16.5$ and $\alpha = {\rm ord}~\mathcal{S}_H < 17$.
\end{Lem}

In \cite{PonSai5}, it was remarked that the coefficient conjecture of Clunie and Sheil-Small is true
if $\mathcal{S}^0_H(\mathcal{S})=\mathcal{S}^0_H$ holds, where
\be\label{eq-SS1}
\mathcal{S}^0_H(\mathcal{S}) = \left\{h+\overline{g} \in \mathcal{S}^0_H :\, h+e^{i \theta}g \in \mathcal{S}~
 \mbox{for some}~~\theta \in \mathbb{R} \right\}
\ee
and $\mathcal{S}^0_H(\mathcal{S})$ contains the class of harmonic mappings convex in one direction. However, this conjecture remains \emph{open}.

Let $\mathscr{F}_H$ be a family of sense-preserving harmonic mappings $f$ with the power series representation as in \eqref{PSSerRep}. Then, the family $\mathscr{F}_H$ is called a linear invariant family, if for each $f\in \mathscr{F}_H$, the function $F$ defined by
$$ F(z) = \frac{f(e^{i\theta}\frac{z+a}{1+\overline{a}z})-f(ae^{i\theta})}{(1-|a|^2)h'(ae^{i\theta})e^{i\theta}}
$$
also belongs to the class $\mathscr{F}_H$ for all $\theta \in \IR$ and $a \in \ID$. A family $\mathscr{F}_H$ is called an affine invariant family, if, in addition, for each $f \in \mathscr{F}_H$, the
function $A(f(z))$ defined by
$$ A(f(z)) = \frac{f(z)+\epsilon \overline{f(z)}}{1+\epsilon \overline{f_{\overline{z}}(0)}}
$$
also belongs to the class $\mathscr{F}_H$ for all $\epsilon \in \ID$.  The order of an affine and linear invariant family $\mathscr{F}_H$ is defined as ${\rm ord}~ \mathscr{F}_H = \sup_{f \in \mathscr{F}_H}|a_2|$. Three well-known affine and linear invariant families are the class $\mathcal{S}_H$, its subclasses $\mathcal{K}_H$ of convex and $\mathcal{C}_H$ of close-to-convex harmonic mappings. It is well known that ${\rm ord}~\mathcal{K}_H = 2$ and ${\rm ord}~\mathcal{C}_H = 3$. In 2004, Starkov \cite{Starkov_2004} (for details see \cite{Starkov_2011}) introduced the order of a linear invariant family (which is not necessarily affine invariant family) $\mathscr{F}_H$ which is defined as follows:
$$ \overline{{\rm ord}}~\mathscr{F}_H = \sup_{f\in \mathscr{F}_H} \frac{|a_2-\overline{b_1}b_2|}{1-|b_1|^2}.
$$
Corresponding to a linear invariant family $\mathscr{F}_H$, we define the family $\mathscr{F}^0_H$ as
$$\mathscr{F}^0_H = \left\{F=\frac{f+\epsilon \overline{f}}{1+\epsilon \overline{f_{\overline{z}}(0)}}:\, f\in \mathscr{F}_H, \epsilon \in \ID, F_{\overline{z}}(0)=0  \right\}.
$$
The following lemma is useful in determining the $\overline{{\rm ord}}~\mathscr{F}_H$.

\begin{Lem}\label{new_ord}{\rm \cite{Starkov_Ganenkova}}
Let $\mathscr{F}_H$ be a linear invariant family of harmonic mappings. Then
$$\overline{{\rm ord}}~\mathscr{F}_H = \sup_{f\in \mathscr{F}_H^0} |a_2| = {\rm ord}~\mathscr{F}^0_H.$$
\end{Lem}

The family $\mathscr{U}_H(\alpha)$ is defined as the union of all affine and linear invariant families $\mathscr{F}_H$ of harmonic functions such that
$\overline{{\rm ord}}~{\mathscr{F}_H} \leq \alpha.$ Set $\mathscr{U}^0_H(\alpha) :=\{f \in \mathscr{U}_H(\alpha): f_{\overline z}(0)=0 \}$. It is now appropriate to state our first main result.

\bthm\label{PS8Thm1}
Let $f= h+\overline{g} \in \mathscr{U}^0_H(16.5)$ with series representation as in \eqref{PSSerRep}. Then we have
\be\label{coeff_bound_all_1}
|a_n| < 5.24 \times 10^{-6} n^{17} ~\mbox{ and }~ |b_n| < 2.32 \times 10^{-7}n^{17} ~~\mbox{for all}~~ n \geq 3.
\ee
\ethm

\br
{\rm
We remark that ${\mathcal S}^0_{H} \subset \mathscr{U}^0_H(16.5)$.
The new bounds in \eqref{coeff_bound_all_1} clearly improves the earlier bounds in \eqref{coeff_bound_starkov}. From the proof of Theorem \ref{PS8Thm1}, we observe that the number $5.24$ in \eqref{coeff_bound_all_1} could be replaced by $4.1006$ for $n \geq 19$ and $2.32$ by $2.25$ for $n \geq 18$. The proof of Theorem \ref{PS8Thm1} relies on the bound $|a_2| \leq 16.5$ for $f\in \mathscr{U}^0_H(16.5)$. If we use the conjectured bound $|a_2| \leq 5/2$, then  Theorem \ref{PS8Thm1} takes an improved version which is stated in Section  \ref{PSS8Sec4}.
}
\er

\subsection{Injectivity of sections of univalent harmonic functions}\label{PSS8SubSec2}
For an analytic function $h(z)=\sum _{k=1}^{\infty}a_k z^k$ in the unit disk $\mathbb{D}$,
the $n$-th section/partial sum $s_n(h)$ of $h$ is defined by
\be\label{AnPSerRep}
s_n(h)(z)=\sum _{k=1}^{n}a_k z^k.
\ee
In \cite{Szego}, Szeg\"{o} proved that every
section $s_n(h)$ of $h\in {\mathcal S}$ is univalent in $|z| < 1/4$ for all $n \geq 2$. The constant $1/4$ is
sharp. If $h \in \mathcal{K}$ (resp. $\mathcal{S}^*$,
and $\mathcal{C}$), then the  $n$-th section $s_n(h)$ is known to be univalent and convex (resp. starlike and close-to-convex)
in the disk $|z| < 1-3n^{-1}\log n$ for all $n \geq 5$ (cf. \cite[Exercise 7, p.~272]{Duren}).
However, the exact radius of univalence $r_n$ of $s_n(h)$, $h \in \mathcal{S}$, remains an \emph{open} problem. By making use of Goluzin's inequality,
Jenkins \cite{Jenkins} proved that $s_n(h)$ is univalent in $|z|<r_n$ for $h\in\mathcal{S}$, where $r_n$
is at least $ 1-n^{-1}(4\log n - \log(4\log n))$ for $n \geq 8$. It is worth pointing out that the result of Jenkins could be improved
if we use de Branges \cite{de_Branges} coefficient estimates $|h^{(n)}(0)/n!|\leq n$ for $h\in {\mathcal S}$.
More precisely, we can easily obtain that $s_n(h)$ is univalent in $|z|<r_n$ for $h\in\mathcal{S}$, where $r_n$ is at least $ 1-n^{-1}(4\log n - 2\log(\log n))$ for $n \geq 7$,
which seems to be the best known radius so far. We avoid the technical details of this fact for obvious reasons.
For related investigations on this topic, see the recent articles \cite{ObSamy13,Hiroshi-Samy-2014, PonSaiStarkov1} and the references therein.

For $f=h+\overline{g} \in \mathcal{S}^0_H$, $n\geq 1$ and $m\geq 2$, the sections/partial sums $s_{n,m}(f)$ of $f$ are defined as
$$s_{n,m}(f)(z)=s_n(h)(z)+\overline{s_m(g)(z)}.
$$
However, the special case $m=n\geq 2$ seems to be interesting in its own merit. In 2013, Li and Ponnusamy \cite{LiSamyNA1, LiSamyNA2, LiSamyCzM} determined the radius of univalence of sections of functions from certain classes of univalent harmonic mappings. For $f$ belonging to ${\mathcal S}_H^{*0}$,  ${\mathcal C}_H^0$, $\mathcal{S}^0_H(\mathcal{S})$ or the class of harmonic mappings convex in one direction, in \cite{PonSaiStarkov1}, the present authors proved that $s_{n,m}(f)$ is univalent in the disk $|z|<r_{n,m}$, where  $r_{n,m}$ is the zero of a rational function.
In the special case $m=n$, $s_{n,n}(f)$ is univalent in the disk $|z|<r_{n,n}$, where
$$r_{n,n} > r^L_{n,n}:=1- \frac{(7\log n  - 4\log(\log n))}{n} ~\mbox{ for }~ n \geq 15.
$$
Moreover, it was also pointed out that $r_{n,m} \geq r^L_{l,l}$,  where  $l=\min\{n,m\}\geq 15$.

In \cite{PonSaiStarkov1}, it was also proved that for $f \in \mathcal{K}^0_H$, each partial sum $s_{n,m}(f)$ is univalent in the disk $|z|<r_{n,m}$, where
$$r_{n,m} \geq 1-\frac{4\log l - 2\log(\log l)}{l} ~ \mbox{ and } ~ l=\min\{n,m\} \geq 7 .
$$
In view of the lack of information on the coefficients of the analytic and co-analytic parts of functions in $\mathcal{S}^0_H$, in contrast to the analytic case, determining the radius of univalence of sections of functions in the class $\mathcal{S}^0_H$ seems to be a difficult task. Nevertheless, in the present article we attempt to consider this problem for the class $\mathcal{S}^0_H$. This is achieved as an application of Theorem \ref{PS8Thm1}.

\bthm\label{PS8Thm2}
Suppose that $f=h+\overline{g} \in \mathcal{S}^0_H$ with the series representation as in \eqref{PSSerRep}.
For $r \in (0, 1)$, define
$$U(r)=\frac{1}{\log r}\left\{-28.5 + \log(r (\log (1/r))^{19})+\log\left[\left(\frac{1-r}{1+r}\right)^{17}
- \left(\frac{1-r}{1+r}\right)^{51}\right] \right\}.
$$
Then, for $n \geq 2$, each section $s_{n,n}(f)$ is univalent in the disk $|z|<r_{n,n}$, where
$$r_{n,n} = \max \{r \in (0, 1):\, 18 \log n = -(n - U(r)) \log r\}.
$$
On the other hand, for fixed $r \in (0.016155, 1)$, $s_{n,n}(f)$ is univalent in $|z| < r$ for all $n \geq N(r)$, where
$$N(r):=\min \{n \ge U(r):\, 18 \log n \leq -(n - U(r)) \log r\}.
$$
\ethm

For example, a routine computation gives the following and so we omit the details.

\bcor\label{PSS8_cor1}
For $f \in \mathcal{S}^0_H$, we have
\begin{enumerate}
\item $s_{n,n}(f)$ is univalent in the disk $|z|<1/4$ whenever $n \geq 81$.
\item $s_{n,n}(f)$ is univalent in the disk $|z|<1/e\approx 0.36788$ whenever $n \geq 131$.
\item $s_{n,n}(f)$ is univalent in the disk $|z|<1/2$ whenever $n \geq 220$.
\end{enumerate}
\ecor

From the proof of Theorem \ref{PS8Thm2}, it is also clear that the result could be improved, if we knew the exact upper bounds on
$|a_n|$ and $|b_n|$ for $f \in \mathcal{S}^0_H$. Therefore it is natural to state an improved form of this result with the assumption
on the order of the family considered. This is done in Section \ref{PSS8Sec4}.

\section{Basic lemmas}\label{PSS8Sec2}

The following results together with Lemma \ref{samy_abu_bound} are useful in the proofs of our main results.

\begin{Lem}{\rm\cite{Starkov}}\label{uni_nec_suf}
A sense-preserving harmonic function $f=h+\overline{g}$ of the form \eqref{PSSerRep} is univalent in $\ID$ if and only if for each
$z \in \ID\setminus\{0\}$ and each $t\in(0, \pi/2]$,
\be\label{PS7inteq2}
\frac {f(re^{i\eta})-f(re^{i\psi})} {re^{i\eta}-re^{i\psi}}
= \sum_{k=1}^{\infty}\left[ (a_k z^k - \overline{b_k z^k})\frac{\sin kt}{\sin t} \right] \ne 0,
\ee
where $a_1=1$, $t=(\eta-\psi)/2$ and $z=re^{i(\eta+\psi)/2}$.
\end{Lem}

\begin{Lem} {\rm\cite{Graf-samy}}\label{two_po_dis}
If $f=h+\overline{g} \in \mathcal{S}^0_H$, $r \in (0, 1)$, $t, \psi \in \IR$, then
\be\label{lower_bound}
\left|\frac {f(re^{it})-f(re^{i\psi})} {re^{it}-re^{i\psi}}\right|\ge \frac {1}{4\alpha r}
\left(\frac {1-r}{1+r}\right)^\alpha \left[1-\left(\frac {1-r}{1+r}\right)^{2\alpha}\right],
\ee
where $\alpha = {\rm ord}~\mathcal{S}_H$.
\end{Lem}

\begin{Lem}\label{Thm_abs_h_g}{\rm \cite{Graf}}
Suppose that $f = h+\overline{g} \in \mathscr{U}_H(\alpha_0)$ with $b_1=f_{\overline{z}}(0)$. For $z \in \mathbb{D}$ with $|z|=r$, $h$ and $g$ satisfy the bounds
$$
|h'(z)| \leq (1+r|b_1|)\frac{(1+r)^{\alpha_0 - 3/2}}{(1-r)^{\alpha_0 + 3/2}} ~\mbox{ and }~
|g'(z)| \leq (r+|b_1|)\frac{(1+r)^{\alpha_0 - 3/2}}{(1-r)^{\alpha_0 + 3/2}}.
$$

\end{Lem}

\section{Proofs of Main Theorems}\label{PSS8Sec3}


\subsection{Proof of Theorem \ref{PS8Thm1}}
Let $f= h+\overline{g} \in \mathscr{U}^0_H(16.5)$. From the power series representation of $h(z)$ given by \eqref{PSSerRep} and Lemma \Ref{Thm_abs_h_g}, we obtain that
\be\label{growth_an}
|a_n| = \left|\frac 1{2\pi i}\int_{|z|=r} \frac {h'(z)} {n z^{n}} dz\right| \leq  \frac {1}{nr^{n-1}} \frac{(1+r)^{15}}{(1-r)^{18}} =: \psi_n(r),
\ee
where $0 < r < 1$. In particular,
$$|a_n| \leq \min_{r \in (0, 1)} \psi_n(r).
$$
In order to obtain the minimum value of the right hand side of the inequality, we need to find the point of minimum of the function
$\log \psi_n(r)$. We see that
$$(\log \psi_n(r))'= \frac{15}{1+r} - \frac{n-1}{r} + \frac{18}{1-r} = 0 \Longleftrightarrow r^2+ \frac {33r}{n+2}-\frac{n-1}{n+2}=0.
$$
It follows that
$$\tau_n = \frac{-33+\sqrt{4n^2+4n+1081}}{2(n+2)}
$$
is the point of minimum and thus,
\be\label{PSS8_eq4}
|a_n|\leq \psi_n (\tau_n) = A(\tau_n) B(\tau_n) ~\mbox{ for all }~ n \geq 3,
\ee
where
$$A(\tau_n)= \left(\frac{2(n+2)}{\sqrt{4n^2+4n+1081}-33}\right)^{n-1} < \left(\frac{2(n+2)}{\sqrt{4n^2+4n+1081}-33}\right)^n
$$
and
$$B(\tau_n)= \frac 1 n \left(\frac{2n + \sqrt{4n^2+4n+1081}-29}{2(n+2)}\right)^{15} \left(\frac{2(n+2)}{2n+37-\sqrt{4n^2+4n+1081}}\right)^{18}.
$$
First, we shall prove that $A(\tau_n) \leq e^{18}$ for all $n \geq 2$.
%
Now, we let
$$ \Psi(x) = \frac{2(x+2)e^{-18/x}}{\sqrt{4x^2+4x+1081}-33}.
$$
Differentiating $\Psi$ with respect to $x$ we get that
$$\Psi'(x) =\frac{6 e^{-18/x}q(x)}{ x^2 t(x)\left(t(x)-33\right)^2},
$$
where $t(x)=\sqrt{4x^2+4x+1081}$ and  $q(x)=q_1(x) - q_2(x)$ with
$$q_1(x) = 12972 + 6534 x + 431 x^2 + 22 x^3 ~\mbox{ and }~ q_2(x) = 396 t(x) + 198 x t(x) + 11 x^2 t(x).
$$
As $q_1(x) > 0 $ and $q_2(x) > 0 $ for $x \ge 0$, it is clear that  $q_1(x) - q_2(x)$ and $q^2_1(x) - q^2_2(x)$ will have the same sign whenever $x \ge 0$.
Computation shows that
$$ q^2_1(x) - q^2_2(x) = 24 (x+2)^2 (44 x^3 + 5381 x^2 + 6438 x -12972) > 0 ~\mbox{ for all }~ x \ge 2,
$$
These observations show that $\Psi'(x) > 0$ for  $x \ge 2$ and hence, $\Psi(x)$ is a increasing function of $x$, whenever $x \ge 2$.
As $\lim_{x \rightarrow \infty} \Psi(x) = 1$, we deduce that $\Psi(x) \leq 1$ for all $x \geq 2$,
which is equivalent to
$$\left(\frac{2(x+2)}{\sqrt{4x^2+4x+1081}-33}\right)^x \leq e^{18} ~\mbox{ for all }~ x \geq 2.
$$
In particular, this observation gives
\be\label{bound_A_tau}
A(\tau_n) \leq e^{18} ~\mbox{  for all $n \geq 2$.}
\ee

Now, we set $n = N+3$, $p(N) = (2N + \sqrt{4N^2+28N+1129} +43)^{18}$ and
$$ T(N) = \sqrt{1+\frac{7}{N}+\frac{1129}{4N^2}}.
$$
A simple calculation shows that
\vspace{6pt}
\noindent
$\ds \frac{B(\tau_{N+3})p(N)}{p(N)}$
\beqq
&=& \frac{2^{-21}3^{-36}N^{17}}{(1+3/N)} \left(1+T(N)+\frac{43}{2N} \right)^3  \left(\frac{1+T(N)+\frac{1}{2N}(17+10T(N)) +\frac{140}{8N^2}}{1+5/N}\right)^{15}\\
& < & \frac{1}{2^{21}3^{36}} \left(1+T(N)+\frac{43}{2N} \right)^3  \left(1+T(N)+\frac{1}{2N}(17+10T(N)) +\frac{140}{8N^2}\right)^{15}N^{17} \\
&\leq& \frac{1}{2^{21}3^{36}} \left(1+T(16)+\frac{43}{32} \right)^3  \left(1+T(16)+\frac{1}{32}(17+10T(16)) +\frac{140}{8 \times {16}^2}\right)^{15}N^{17}
\eeqq
for all $N \geq 16$. Since $T(16) \approx 1.59375$, the last inequality then gives that
\beq\label{bound_B_tau}
B(\tau_{N+3}) &\leq& (3.17691 \times 10^{-24}) \times 61.0466 \times (3.22016 \times 10^8) \nonumber \\
&\leq& 6.2452 \times 10^{-14}  N^{17} \mbox{ for all $ N \geq 16$.}
\eeq
Hence, by \eqref{bound_A_tau} and \eqref{bound_B_tau}, one obtains that
$$A(\tau_{N+3})B(\tau_{N+3}) \leq e^{18} \times 6.2452 \times 10^{-14}  N^{17} \approx 4.1006 \times 10^{-6} N^{17} ~\mbox{ for all $ N \geq 16$}.
$$
By a direct but lengthy computation or by Mathematica, we can easily see that
$$A(\tau_n) B(\tau_n) \leq  5.24 \times10^{-6} n^{17} ~\mbox{ for }~ 3 \leq n \leq 18.
$$
Therefore, using these two estimates, the inequality \eqref{PSS8_eq4} reduces to
$$|a_n| \leq 5.24 \times10^{-6} n^{17} ~\mbox{  for all $n \geq 3$.}
$$
Similarly, from the power series representation of $g$ given by \eqref{PSSerRep}, one sees that
$$|b_n| = \left|\frac 1{2\pi i}\int_{|z|=r} \frac {g'(z)} {n z^{n}} dz\right| \leq  \frac {1}{nr^{n-2}} \frac{(1+r)^{15}}{(1-r)^{18}} : =\phi_n(r) = \frac{n-1}{n}\psi_{n-1}(r),
$$
where $0 < r < 1$ and $\psi_n(r)$ is defined as in \eqref{growth_an}. In particular,
$$|b_n| \leq \min_{r \in (0, 1)} \phi_n(r).
$$
Using similar arguments as above, we get that
$$|b_n| \leq A_1(\rho_n) B_1(\rho_n) ~\mbox{ for all }~ n \geq 3,
$$
where
$$ \rho_n =  \frac{-33+\sqrt{4n^2-4n+1081}}{2(n+1)},
$$
$$A_1(\rho_n)= \left(\frac{2(n+1)}{\sqrt{4n^2-4n+1081}-33}\right)^{n-2} = A(\tau_{n-1}) \leq e^{18} ~\mbox{ for all }~ n \geq 3
$$
and
$$B_1(\rho_n)= \frac 1 n \left(\frac{2n + \sqrt{4n^2-4n+1081}-31}{2(n+1)}\right)^{15} \left(\frac{2(n+1)}{2n+35-\sqrt{4n^2-4n+1081}}\right)^{18}.
$$
Setting
$$ l(n) = \sqrt{1+\frac{1081}{4n^2}-\frac 1 n},
$$
we obtain that
$$B_1(\rho_n)= \frac{8}{n 144^{18}(n+1)^{15}} \left[8n^2(1+l(n))+4n(1+2l(n)) - 4\right]^{15} \left[2n(1+l(n))+35 \right]^3.
$$
Using the fact that $\sqrt{1+x} \leq 1 + \sqrt{x}$ for $x \geq 0$, we get that
$$ B_1(\rho_n)= \frac{n^{17}}{2^{21}3^{36}} \left(2 + \frac{18}{n} + \frac{16}{n^2} \right)^{15} \left(2 + \frac{34}{n} \right)^3 \leq 3.425 \times 10^{-15} n^{17}
$$
for all $n \geq 18$. Therefore,
$$A_1(\rho_n)B_1(\rho_n) \leq  3.425 \times 10^{-15} e^{18} n^{17} \approx 2.25 \times 10^{-7}  n^{17} ~\mbox{ for all } n \geq 18.
$$
By a direct computation with the help of Mathematica, one can see that
$$A_1(\rho_n)B_1(\rho_n) \leq 2.32 \times 10^{-7}  n^{17} ~\mbox{ for }~ 3 \leq n \leq 18.$$
Therefore, $|b_n| \leq 2.32 \times 10^{-7}  n^{17} $ for all $n \geq 3$.
\hfill $\Box$

\subsection{Proof of Theorem \ref{PS8Thm2}}
Suppose that $f=h+\overline{g}$ belongs to ${\mathcal S}_H^{0}$. Set $F_r(z)=f(rz)/r$ for $0 < r < 1$.
Then  $F_r(z) \in \mathcal{S}^0_H$. In view of Lemma \Ref{uni_nec_suf}, it is clear that $s_{n,m}(f)$
is univalent in $|z|<r$ if and only if $s_{n,m}(F_r)(z)$ is sense-preserving in $\ID$ and the associated harmonic polynomial $P_{n,m,r}(z)$
has the property that
$$P_{n,m,r}(z) := \sum_{k=1}^{\infty}\left[ (a'_k z^k - \overline{b'_k z^k})\frac{\sin kt}{\sin t} \right] \ne 0
~\mbox{ for all } z\in\ID\setminus\{0\} \mbox{ and } t\in(0,\pi/2],
$$
where
$$ a'_k = a_k r^{k-1} ~\mbox{ for }~ k\in \{1,2,\ldots,n\} ~\mbox{ and }~ a'_k = 0 ~\mbox{ if }~ k >n
$$
and
$$ b'_k = b_k r^{k-1} ~\mbox{ for }~ k\in \{1,2,\ldots,m\} ~\mbox{ and }~ b'_k = 0 ~\mbox{ if }~ k >m.
$$
Now, we set $t=(\eta-\psi)/2$ and $z=\rho e^{i(\eta+\psi)/2}\in\ID$ in \eqref{PS7inteq2}. Note that the function in the right side of the inequality \eqref{lower_bound} in Lemma \Ref{two_po_dis} decreases with increasing value of $\alpha$, where $\alpha = {\rm ord}\,\mathcal{S}_H$. As $F_r \in \mathcal{S}^0_H$ and $\alpha < 17$, we apply Lemma \Ref{two_po_dis}
to the function $F_r$ and get that
$$ 
\left|\sum_{k=1}^{\infty}\left[ (a_k z^k - \overline{b_k z^k})r^{k-1}\frac{\sin kt}{\sin t} \right]\right| \geq  \frac {1}{68 r}
\left(\frac {1-r}{1+r}\right)^{17} \left[1-\left(\frac {1-r}{1+r}\right)^{34}\right].
$$ 
In order to find a lower bound for $|P_{n,m,r}(z)|$, we need to find an upper bound for
$$\left|R_{n,m,r}(z)\right| = \left| \sum_{k=n+1}^{\infty}\left[ a_k r^{k-1} z^k \frac{\sin kt}{\sin t} \right] - \sum_{k=m+1}^{\infty}\left[\overline{(b_k r^{k-1} z^k)} \frac{\sin kt}{\sin t} \right]\right|.
$$
Using Theorem \ref{PS8Thm1} and the fact that $|\sin kt| \le k \sin t$  for all $t \in [0, \pi/2]$ and $k \in \IN$, we get that
\beq\label{PSS8_eq2}
|R_{n,m,r}(z)| &\leq& \sum_{k=n+1}^{\infty}5.24 \times 10^{-6} k^{18} r^{k-1} + \sum_{k=m+1}^{\infty} 2.32 \times 10^{-7} k^{18} r^{k-1}  \\ \nonumber
 &=:&R_{n, r} + T_{m, r}.
\eeq
Set $\psi(n,m,r) = C_{17}(r) - (R_{n, r} + T_{m, r})$, where
$$ C_{17}(r) = \frac {1}{68 r} \left(\frac {1-r}{1+r}\right)^{17} \left[1-\left(\frac {1-r}{1+r}\right)^{34}\right].
$$
The inequality $|P_{n,m,r}(z)|>0$ holds for all $z\in \ID\setminus\{0\}$, whenever $\psi(n,m,r)>0$.
In \cite[Lemma 1]{Graf-samy} it was shown that $r \mapsto C_{\alpha}(r)$ is strictly decreasing on $(0, 1)$. This fact implies that
$\psi(n,m,r)$ is decreasing in $(0, 1)$ and thus $\psi(n,m,r)>0$ for all $r\in(0, r_{n,m})$, where $r_{n,m}$ is the unique positive root of the equation $\psi(n,m,r)=0$ which is less than $1$. It is easy to see that $s_{n,m}(F_r)(z)$ is sense-preserving in $\ID$ provided $r\in(0, r_{n,m})$ (see e.g \cite{PonSaiStarkov1}) and hence, $s_{n,m}(f)$ is univalent in $|z| < r_{n,m}$.

Now, let us consider the special case $m=n$. In this case $\psi(n,m,r)$ reduces to
$$\psi(n,n,r)=\frac {1}{68 r} \left(\frac {1-r}{1+r}\right)^{17} \left[1-\left(\frac {1-r}{1+r}\right)^{34}\right]
- \sum_{k=n+1}^{\infty} 54.72 \times 10^{-7} k^{18} r^{k-1}.
$$
From our discussion, it is clear that $s_{n,n}(f)$ is univalent in $|z|<r_{n,n}$, where $r_{n,n}$ is the unique positive root of the equation $\psi(n,n,r)=0$ which is less than $1$. In order to compute the lower bound for $r_{n,n}$, we consider the function $T(x) = x^{18}r^{x-1}$. It follows that $T(x)$ is a decreasing function of $x$ in the interval $[-18/\log r, \infty)$.
Whenever $-n \log r > 18$, we have
$$ \sum_{k=n+1}^{\infty} k^{18} r^{k-1} < \int_{n}^{\infty} T(x)\, \mathrm{d}x = \int_{n}^{\infty} x^{18}r^{x-1}\, \mathrm{d}x.
$$
Applying integration by parts repeatedly, we obtain that
$$\int_{n}^{\infty} x^{18}r^{x-1}\, \mathrm{d}x = \frac{1}{r}\left\{\frac{n^{18}r^{n}}{|\log r|} + \frac{18n^{17}r^{n}}{|\log r|^2} +  \frac{18 \times 17n^{16}r^{n}}{|\log r|^3} + \cdots +  \frac{18!r^{n}}{|\log r|^{19}} \right\}.
$$
Choose $n$ large enough so that $-n \log r > 18a \log n \geq 18$, where $a \in (1, \infty)$, which implies that
$$ r^{n} n^j \leq r^{n} n^{18} \leq r^{n((a-1)/a)} ~\mbox{ for all }~ j=1, 2, \ldots, 18.
$$
Hence, we get that
\beqq
\int_{n}^{\infty} x^{18}r^{x-1}\, \mathrm{d}x &\leq& \frac{18! ~r^{n((a-1)/a)-1}}{|\log r|^{19}} \left\{\frac{|\log r|^{18}}{18!} + \frac{|\log r|^{17}}{17!} +  \cdots +  \frac{|\log r|}{1!} + r^{n/a} \right\}\\
&\leq& \frac{18! ~r^{n((a-1)/a)-2}}{|\log r|^{19}}.
\eeqq
Using the above inequality and \eqref{PSS8_eq2}, we obtain that
\beq\label{PSS8_eq3}
|R_{n,n,r}| &\leq& R_n + T_n = \sum_{k=n+1}^{\infty}(52.4+2.32) 10^{-7} k^{18} r^{k-1}\\ \nonumber
&\leq&  54.72 \times 10^{-7} ~\frac{18! ~r^{n((a-1)/a) - 2} }{|\log r|^{19}}
\eeq
provided $\ds r \leq n^{-18a/n}$ for some $a \in (1, \infty)$.
Therefore, $\psi(n,n,r) > 0$ whenever
$$\frac {1}{68} \left(\frac {1-r}{1+r}\right)^{17} \left[1-\left(\frac {1-r}{1+r}\right)^{34}\right]
- 54.72 \times 10^{-7} \frac{r^{n((a-1)/a) - 1}18!}{|\log r|^{19}} \geq 0.
$$
This gives that
$$ u(a, r):= \frac{a}{a-1} U(r) \le n,
$$
where
$$U(r)=\frac{1}{\log r}\left\{-28.5 + \log(r |\log r|^{19})+\log\left[\left(\frac{1-r}{1+r}\right)^{17}- \left(\frac{1-r}{1+r}\right)^{51}\right] \right\}.
$$
The lower bound for $r_{n,n}$ is obtainable for all $n \geq 2$ and this follows from the fact that $aU(r)/(a-1) \rightarrow 2^-$ as $a \rightarrow \infty$ and $r \rightarrow 0^+$. Similarly, $U(r) \rightarrow \infty$ as $r \rightarrow 1^-$. Therefore, $U(r)$ accepts all values from $(2-\delta, \infty)$ if $r \in (0, 1)$, where $\delta$ is some positive constant.

%

From the above discussion, it is clear that $s_{n,n}(f)$ is univalent in $|z|<r$, whenever $u(a, r)\leq n$ and $r \leq n^{-18a/n}$ for some $a \in (1, \infty)$. The inequality $u(a, r)\leq n$ holds for any $r \in (0, 1)$, such that $U(r) < n$, if we choose $a = a^* = n/(n-U(r))$. The inequality $r \leq n^{-18a/n}$ holds true if we choose $r = R_{n,n}$ with $a = a^*$, where
$$
R_{n,n} = \max \{r \in (0, 1):\, 18 \log n \le -(n - U(r)) \log r ~\}.
$$
In this expression, the maximum is reached, because for fixed $n$, $-(n - U(r)) \log r \rightarrow -\infty$ as  $r \rightarrow 1^-$ and $-(n - U(r)) \log r \rightarrow \infty$ as $r \rightarrow 0^+$. For every $n \ge 2$, the function continuously depends on $r$. Thus,
$$r_{n,n} = \max \{r \in (0, 1):\, 18 \log n = -(n - U(r)) \log r ~\}.
$$

We remark here that $r \mapsto U(r)$ is strictly increasing on $(0.016155, 1)$. In order to prove that, it is enough to show that
$$U_1(r)= - \log(r |\log r|^{19})-\log\left[\left(\frac{1-r}{1+r}\right)^{17}- \left(\frac{1-r}{1+r}\right)^{51}\right]
$$
is strictly increasing on $(0.016155, 1)$. A computation shows that
$$ U'_1(r) = \frac{1}{r}\left(-1 - \frac{19}{\log r} \right) + \frac{34}{1-r^2}\left( 1 - \frac{2}{\left(\frac{1+r}{1-r}\right)^{34} - 1} \right) > 0 ~\mbox{for}~ r \ge 0.016155.
$$



Next, for a given $r \in (0.016155, 1)$, we consider the problem of finding the least positive integer $N(r)$ such that $s_{n,n}(f)$ is univalent in the disk $|z|<r$ for all $n \geq N(r)$. In order to guarantee the univalency of $s_{n,n}(f)$ for all $n \geq N(r)$, the number $N(r)$ must be greater than or equal to $u(a, r)$ and $r \leq n^{-18a/n}$ for all $n \ge N(r)$. Using the above arguments, we obtain that
$$ N(r) = \min\{n \ge U(r):\, 18 \log n \leq -(n - U(r)) \log r \}.
$$
This completes the proof. \hfill $\Box$
%

\section{Concluding remarks}\label{PSS8Sec4}
It is known that the inequality $|a_2| \leq 5/2$ holds for functions in $\mathcal{C}^0_{H}$ and from Lemma \Ref{new_ord}, it follows that $\mathcal{C}_{H} \subset \mathscr{U}_H(5/2)$. Theorems \ref{spec_coeff_bound} and \ref{coeff_bound_all_1} below are the analog of Theorem \ref{PS8Thm1} and \ref{PS8Thm2} for the families $\mathscr{U}^0_H(5/2)$ and $\mathscr{U}^0_H(5/2) \cap \mathcal{S}^0_H$, respectively. If $|f''(0)/2| \leq 5/2$ for all harmonic mappings $f \in \mathcal{S}^0_{H} $ as conjectured by Clunie and Sheil-Small, then $\mathcal{S}^0_{H} \subset \mathscr{U}^0_H(5/2)$.

\bthm\label{spec_coeff_bound}
Suppose that $f= h+\overline{g} \in \mathscr{U}^0_H(5/2)$ with series representation as in \eqref{PSSerRep}. Then for all $n \geq 3$,
$$|a_n| \leq \frac{8 (2 + n)^3 (\sqrt{4 n^2+ 4 n+17}+ 2 n-1)}{n (\sqrt{4 n^2+ 4 n+17}- 2 n-9)^4} \left(\frac{\sqrt{4 n^2+ 4 n+17}-5}{2(2+n)}\right)^{1-n}
$$
and
$$|b_n| \leq \frac{8 (1 + n)^3 (\sqrt{4 n^2- 4 n+17}+ 2 n-3)}{n (\sqrt{4 n^2- 4 n+17}- 2 n-7)^4} \left(\frac{\sqrt{4 n^2- 4 n+17}-5}{2(1+n)}\right)^{2-n}.
$$
In particular, the following bounds hold:
$$ 
|a_n| \leq \frac{3n^3}{4}  ~\mbox{ and }~ |b_n| \leq \frac{43 n^3}{100} ~\mbox{ for all }~ n \geq 3.
$$ 
\ethm

As the proof of Theorem \ref{spec_coeff_bound} is similar to that of the proof of Theorem \ref{PS8Thm1}, we omit the details here. As an application of Theorem \ref{spec_coeff_bound}, we prove the following result:

\bthm\label{part_spec_case}
Suppose that $f=h+\overline{g} \in \mathscr{U}^0_H(5/2) \cap \mathcal{S}^0_H$. Then the partial sums $s_{n,m}(f)$ is univalent in the disk $|z|<r_{n,m}$. Here $r_{n,m}$ is the unique positive root of the equation $\varphi(n,m,r)=0$, where
\be\label{PSS8_thm3_eq1}
\varphi(n,m,r) = \frac {1}{12 r}
\left(\frac {1-r}{1+r}\right)^3 \left[1-\left(\frac {1-r}{1+r}\right)^{6}\right] - R_{n,r} - T_{m,r},
\ee
with
$$ 
R_{n,r} = \sum_{k=n+1}^{\infty}\frac{3k^4}{4} r^{k-1} ~\mbox{ and }~ T_{m,r} = \sum_{k=m+1}^{\infty}\frac{43k^4}{100} r^{k-1}.
$$
In particular, each section $s_{n,n}(f)$ is univalent in the disk $|z|<r_{n,n}$, where
\be\label{PSS8_eq5}
r_{n,n} > r^L_{n,n}:=1- \frac{(8\log n  - 4\log(\log n))}{n} ~\mbox{ for }~ n \geq 20.
\ee
Moreover, $r_{n,m} \geq r^L_{l,l}$,  where  $l=\min\{n,m\}\geq 20$.
\ethm
\bpf
The first part of the proof is similar to the proof of Theorem \ref{PS8Thm2}. Following the proof of Theorem \ref{PS8Thm2}, under the hypothesis of Theorem \ref{part_spec_case}, we get that
$s_{n,m}(f)$ is univalent in the disk $|z|<r_{n,m}$, where $r_{n,m}$ is the unique positive root of the equation $\varphi(n,m,r)=0$.
Here $\varphi(n,m,r)$  is given by \eqref{PSS8_thm3_eq1}.
For $m=n$, we have
$$ 
\varphi(n,n,r)=\phi(r)  - (R_{n,r} + T_{n,r}),
$$ 
where
$$ \phi(r) = \frac{(1-r)^3 (3 + 10 r^2 + 3 r^4)}{3(1+r)^9}
$$
and
\beqq
 R_{n,r} + T_{n,r} &=& \sum_{k=n+1}^{\infty} \frac{59}{50}k^{4} r^{k-1}\\
 & =& \frac{59 r^n}{50(1-r)^5} \left\{1 + 4 n^3 (1 - r)^3 + n^4 (1 - r)^4 + 11 r + 11 r^2 + r^3 \right. \\
 &~~& \left. + 6 n^2 (1 - r)^2 (1 + r) - 4 n (r^3+ 3 r^2- 3 r-1 )\right\}.
\eeqq
The inequality $\varphi(n,n,r) \geq 0$ holds if and only if $0 <  k(n, r) \le 1$, where
$$ k(n, r) := \frac{R_{n,r} + T_{n,r}}{\phi(r)}.$$
Now, we show that for every fixed integer $n \ge 2$, $k(n, r)$ is a increasing function of $r$ in the interval $~[0, 1)$. In order to show that
$k(n, r)$ is a increasing function of $r$, it is enough to show that $\phi_1(r) = (1-r)^3 (3 + 10 r^2 + 3 r^4)$ is decreasing function of $r$ in the interval $[0, 1)$. Since
$$\phi_1'(r) = (1-r)^2 (-9 + 20 r - 50 r^2 + 12 r^3 - 21 r^4) \leq (1-r)^2 (-9 + 20 r - 38 r^2) < 0
$$
for all $r \in [0, 1)$, $\phi_1(r)$ is decreasing and hence $k(n, r)$ is increasing function and $k(n, r) > 0 $ for $r \in (0, 1)$. As $\lim_{n\rightarrow\infty}(R_{n,r} + T_{n,r}) =0$, it is clear that the radius $r_{n,n}$ of univalence approaches $1$. This suggests that $r_{n,n} \geq r_{n,n}^L:=1-x_n/n$, where $x_n$ is positive and increasing sequence of real numbers such that $x_n=o(n)$.

Let us compute the approximate value of $r_{n,n}$ for large values of $n$. By setting $r = 1-x/n$ in $k(n, r)$, and making use of the fact that $(1-x/n)^n \leq e^{-x}$ for $x \geq 0$, we get that $k(n, 1-x/n) \leq t(x, n)$, where
$$t(x, n) := \frac{177e^{-x}n^8}{50 x^8}\left(2 - \frac{x}{n}\right)^9\frac{q(x,n)}{16 n^4 - 32 n^3 x + 28 n^2 x^2 - 12 n x^3 + 3 x^4},
$$
with
$$q(x, n) = n[n^3 (24 + 24 x + 12 x^2 + 4 x^3 + x^4) - 6 n^2 x (6 + 4 x + x^2)+ 2 n x^2 (7 + 2 x) -x^3].
$$
We may set $x_n=8\log n  - 4\log(\log n)$ and we observe that $1 - x_n/n > 0$ only when $n \geq 20$. Therefore, we consider the case $n \ge 20$.
In order to show that $\varphi(n,n,r) \geq 0$ for all $r \in (0, 1-x_n/n)$, it suffices to prove that  $t(x_n, n) \leq 1$ for all $n \ge 20$.

Set
$$ t(x_n, n) = T_1(n) T_2(n) T_3(n),
$$
where
\beqq
T_1(n) &=& \frac{177}{50 \times 2^7 \times \left(2 - \frac{\log(\log n)}{\log n}\right)^8} \left(1 - \frac{4\log n - 2 \log(\log n)}{n} \right)^9\\
       &\le& \frac{177}{50 \times 2^7 \times\left(2 - \frac{\log(\log 20)}{\log 20}\right)^8} \approx 0.000555~\mbox{ for all }~ n \ge 20,
\eeqq
\beqq
T_2(n) &=&  \frac{q(x_n, n)}{(n \log n)^4} \le \frac{24}{(\log n)^4} + \frac{8}{(\log n)^3}\left(24-\frac{36}{n} \right)+ \frac{8^2}{(\log n)^2}\left(\frac{14}{n^2} - \frac{24}{n}\right)\\
       &~~& \hspace{2cm} +~~~~ \frac{8^3}{(\log n)}\left(4 - \frac{6}{n} +\frac{4}{n^2} - \frac{1}{n^3}\right) + 8^4\\
       &\le& \frac{24}{(\log n)^4} + \frac{8\times 24}{(\log n)^3}  + \frac{8^3 \times 4}{(\log n)} + 8^4 = : S_2(n)\\
       &\le& S_2(20) \approx 4787.08 ~\mbox{ for all }~ n \ge 20
\eeqq
and
\beqq
T_3(n) &=&  \frac{1}{16 - 32 (x_n/n) + 28 (x_n/n)^2 - 12 (x_n/n)^3 + 3 (x_n/n)^4} = : S_3(x_n/n)\\
       &\le& S_3(x_{20}/20) \approx 0.333 ~\mbox{ for all }~ n \ge 20.
\eeqq
Therefore,
$$ t(x_n, n) \le 0.000555 \times  4787.08 \times  0.333 \approx 0.885 < 1 ~\mbox{ for all }~ n \ge 20.
$$
This completes the proof of the theorem.
\epf


A rough estimate on $r_{n,n}$ gives the following  which may be compared with Corollary \ref{PSS8_cor1}.

\bcor\label{PSS8_cor2}
Suppose that $f=h+\overline{g} \in \mathscr{U}^0_H(5/2) \cap \mathcal{S}^0_H$. Then $s_{n,n}(f)$ is univalent in the disk
$|z|<r$, where {\rm (i)} $r=1/4$ whenever  $n \geq 10$;  {\rm (ii)}  $r=1/2$ whenever $n \geq 29$; and  {\rm  (iii)} $r=3/4$ whenever $n \geq 98$.
\ecor

Better lower bounds for the radius of univalence $r_{n,n}$ of $s_{n,n}(f)$  (under the assumptions of Theorem \ref{part_spec_case})
for certain values of $n$ are listed in Table \ref{tab2}. They are obtained by
solving the equation $\varphi(n,n,r)=0$.

\begin{table}[htp]
\begin{center}
\begin{tabular}{|l|l||l|l|l|}
  \hline
  Value of $n$ & Lower bound for $r_{n,n}$ &Value of $n$ &Lower bound for $r_{n,n}$\\
  \hline
  2 & 0.0635798 &  10 &  0.269796\\
  \hline
  3 &  0.0952634 & 50 &  0.625779\\
  \hline
  4 &  0.12535 & 100 &  0.753905\\
  \hline
  5 &  0.153603&  354 &  0.900055\\
  \hline
\end{tabular}
\end{center}
\caption{Values of $r_{n,n}$ for certain values of $n$\label{tab2}}
\end{table}

\subsection*{Acknowledgements}
The research of the first author was supported by the project RUS/RFBR/P-163 under Department of Science \& Technology (India) and this
author is currently on leave from Indian Institute of Technology Madras. The work of the second author was supported
by NBHM (DAE), India. The third author is supported by Russian Foundation for Basic Research (project 17-01-00085) and the Strategic 
Development Program of Petrozavodsk State University.

\end{document}